\documentclass{article}

 \UseRawInputEncoding

\newtheorem{theorem}{Theorem} 

\newcommand\tpv{{T_{\mbox{PV}}}}

\newcommand\qed{\begin{flushright} {\bf q.e.d.} \end{flushright} }
\newcommand\prf{\noindent {\bf Proof :}}  

\newcommand\bits{\{0,1\}}
\newcommand\nn{{\bits^n}}

\newcommand\mm{{\bits^m}}
\newcommand\kk{{\bits^k}}

\newcommand\tru{{\mbox{{\bf tt}}}}
\newcommand\sz{{\mbox{Size}}}

\begin{document}

\title{Small circuits and dual weak PHP in the universal theory of p-time algorithms}

\author{Jan Kraj\'{\i}\v{c}ek}

\date{Faculty of Mathematics and Physics\\
Charles University\thanks{
Sokolovsk\' a 83, Prague, 186 75,
The Czech Republic, {\tt krajicek@karlin.mff.cuni.cz}}}

\maketitle

\begin{abstract}
We prove, under a computational complexity hypothesis, that it is consistent with the true universal
theory of p-time algorithms that a specific p-time function extending $n$ bits to $m \geq n^2$ bits
violates the dual weak pigeonhole principle: every string $y \in \mm$ equals the value 
of the function for some $x \in \nn$.
The function is the truth-table function assigning to a circuit the table of the function it computes
and the hypothesis is
that every language in P has circuits of  a fixed polynomial size $n^d$.

\end{abstract}

Consider a first-order language having a function symbol for every deterministic
p-time clocked Turing machine, the symbol
being interpreted over $\mathbf N$ by the function computed by the machine. 
Relations computable by p-time machines are formally represented by their characteristic functions
but we write, for example, $x \le y$ instead of $\le(x,y)=1$.
We shall denote this language $L_{\mbox{PV}}$
and the theory of all true universal sentences in the language by $T_{\mbox{PV}}$; the notation alludes to
the influential theory PV  
introduced by Cook \cite{Coo75} although its language is defined in a much more complicated
way (because of the intended links with proof complexity, cf.\cite{Coo75} or \cite[Chpt.5]{kniha}). Note that
a number of important theorems from computational complexity, including the PCP
theorem or all valid instances of the NP-completeness of SAT 
(where there are p-time functions sending witnesses to witnesses)
or various lower bounds, can be expressed as universal 
statements in the language and hence are axioms of theory $T_{\mbox{PV}}$ (cf. \cite[Sec.22.3]{prf}). 

Dual weak pigeonhole principle (dWPHP) for a function $g$ says that for no $n < m$ can $g$ map $n$-bit strings
onto all $m$-bit strings. It was first considered in the context of bounded arithmetic by Wilkie
who proved a witnessing theorem for a particular theory having dWPHP among its axioms; this was written up in
\cite[Thm.7.3.7]{kniha}. The theory was suggested as the {\em basic theory} (BT) for formalization of complexity theory
(and, in particular, of probabilistic constructions) in \cite{Kra-wphp} and indeed 
Je\v r\'abek \cite{Jer-phd,Jer04,Jer-apc1} succeeded in it spectacularly. The dWPHP is also linked with proof
complexity (cf. \cite{Kra-dual}) and one can view the question (posed in \cite{Kra-wphp})
whether one of the theories related to p-time algorithms proves dWPHP for p-time functions
as a uniform version of particular propositional lengths-of-proofs problems (about proof complexity generators,
cf.\cite[Sec.19.4]{prf}). 
We do not give definitions of the notions mentioned above or details of the statements as they are 
not technically relevant to this paper and serve here only as a motivation for the expert reader. The non-expert
reader can find this background in the references given and, in particular, all of it in \cite{prf}.

The unprovability of dWPHP in $\tpv$ implies that $\mbox{P} \neq \mbox{NP}$: Paris, Wilkie and Woods \cite{PWW} proved
that dWPHP for p-time functions is provable in a theory
having induction for all predicates in the p-time hierarchy (theory $T_2(\mbox{PV})$ of Buss \cite{Bus-book}) but
that theory would follow from $\tpv$ if it were that P = NP. This is because if satisfiability can be solved by p-time 
algorithm $f$ that statement is universal:
$$
Sat(x,y) \rightarrow Sat(x, f(x))
$$
and hence in $\tpv$ (here $Sat(x,y)$ formalizes that $y$ is a satisfying assignment for formula $x$).
But then every bounded formula is provably in $\tpv$ equivalent to an open formula and hence induction for
all bounded formulas follows from induction for open formulas which is in $\tpv$ provable via usual binary search.
This means that if we want to prove that dWPHP is not provable in $\tpv$
we ought to expect to use some hypothesis that itself implies $\mbox{P} \neq \mbox{NP}$. The hypothesis
we shall use is the following:

\medskip
\noindent
{\bf Hypothesis (H):} 

{\em There exists constant $d \geq 1$ such that every language in P can be decided by circuits of size $O(n^d)$:
$\mbox{P}\subseteq \sz(n^d)$.}

\medskip

The popular expert opinion finds (H) unlikely, I suppose, but the reader should note that there are 
no technical results that would support the skepticism. In anything, the notorious unability to prove even
$10n$ lower bound for general circuits suggest that the possibility that (H) is true
cannot be simply dismissed. It has also the attractive feature that it implies $\mbox{P} \neq \mbox{NP}$
(there are languages in the polynomial-time hierarchy that have no size $O(n^d)$ circuits, cf. Kannan's theorem \cite{Kan}).
Thus, in principle, one could prove $\mbox{P} \neq \mbox{NP}$ by proving circuit upper bounds rather than
by proving lower bounds. It is less attractive that it implies also $\mbox{E}\subseteq \sz(2^{o(n)})$ 
(a language in E becomes p-time computable if the inputs are padded)
and hence it disproves the foundational hypothesis of universal derandomization. But this is not an a priori reason to abandon
(H) as $\mbox{E}\not\subseteq \sz(2^{o(n)})$ is itself only a hypothesis. On the other hand 
(H) is good for proof complexity:
together with \cite[Thm.2.1]{Kra-di} the statement $\mbox{E}\subseteq \sz(2^{o(n)})$ (and hence (H)) imply that either
$\mbox{NP}\neq\mbox{coNP}$ or that there is no p-optimal propositional proof system; proving (or disproving) one of these two
statements are the two fundamental problems of proof complexity.
The hypothesis (with linear size circuits) is often attributed to Kolmogorov, see the discussion in \cite[Sec.20.2]{Jukna-book}. 

\medskip

Next we need to define the {\bf truth-table function} $\tru_{s,k}$. It takes as an input 
a circuit with $k$ inputs of size $\le s$ and outputs its truth table, $2^k$ bits. 
A size $\le s$ circuit can be encoded by, say, $10s\log s$ bits exactly 
and hence for $10 s \log s < 2^k$ this is a function from a smaller 
set into a bigger one. Our size function $s = s(k)$ will have the form $s(k) := 2^{\epsilon k}$ for some fixed 
$0 < \epsilon < 1$. Hence any such $\tru_{s,k}$ is a p-time function.

\begin{theorem} \label{8.4.20a}
{\ }

Assume hypothesis (H). Then for every $0 < \epsilon < 1$ and $s = s(k) := 2^{\epsilon k}$
the theory $\tpv$ does not prove the sentence
\begin{equation} \label{9.4.20a}
\forall 1^m (m=2^k > 1) \exists y \in \mm \forall x \in \nn,\ \tru_{s,k}(x) \neq y
\end{equation}
expressing the dWPHP for $\tru_{s,k}$, where $n := 10 s \log s$.
\end{theorem}

\prf

Assume that $\tpv$ proves (\ref{9.4.20a}). By the KPT theorem (cf.\cite{KPT} or \cite[Thm.7.4.1]{kniha} or 
\cite[Cor.12.2.4]{prf}) there a p-time functions
\begin{equation} \label{10.4.20a}
f_1(z), f_2(z,w_1), \dots, f_t(z,w_1, \dots, w_{t-1})
\end{equation}
such that for any $m = 2^k > 1$ and
any $b_1, \dots, b_t, C_1, \dots, C_{t-1}$:

\begin{itemize}

\item either $b_1 \notin rng(\tru_{s,k})$ for $b_1 = f_1(1^m) \in \mm$ or, if $b_1 \in rng(\tru_{s,k})$
and $b_1 = \tru_{s,k}(C_1)$, 

\item $b_2 \notin rng(\tru_{s,k})$ for $b_2 = f_{2}(1^m, C_1) \in \mm$ or, if $b_2 \in rng(\tru_{s,k})$
and $b_2 = \tru_{s,k}(C_2)$,

\item $\dots$, or

\item $b_t \notin rng(\tru_{s,k})$ for $b_t = f_{t}(1^m, C_1, \dots, C_{t-1}) \in \mm$.

\end{itemize}
Define constants $\delta_i := (2d)^{-i}$, for $i = 0, \dots, t$, and parameters $m_i := m^{\epsilon \delta_i}$	
where $d$ is the constant from (H) and $m$ is large enough.

We first show that $f_1$ cannot find a suitable $b_1$. Define the function $\hat f_1$ that has $m_t + k$ variables
and on inputs $1^{m_t}$ and $i \in \kk$ computes the $i$-th bit of $f_1(1^m)$. The string $1^{m_t}$
has the only purpose to make $\hat f_1$ p-time. By hypothesis (H) there is a circuit $C'_1(z,i)$ with the same variables
as $\hat f_1$ that computes $\hat f_1$. Define $C_1$ by substituting $1^{m_t}$ for $z$ in $C'_1$ and leaving just
the $k$ variables for bits of $i$. Note that $C_1$ has size $O((m_t + k)^d)$ and thus can be encoded
by $\le m_{t-1}$ bits. Further, by its definition, $\tru_{s,k}(C_1) = b_1$.

Now we show that $f_2$ does not compute a
suitable $b_2 := f_2(1^m, C_1)$ either. As before define function $\hat f_2$
that now takes three inputs: string $1^{m_{t-1}}$, circuit $C_1$ (substituted for $w_1$) and $i \in \kk$, 
and computes the $i$-th bit of $f_2(1^m, C_1)$.
Applying (H) we get a circuit $C'_2$ with the same $2 m_{t-1} + k$ variables
as $\hat f_2$ that computes the function. Define $C_2$ by substituting $1^{m_{t-1}}$ for $z$ and bits defining $C_1$ for $w_1$
in $C'_2$, and leaving just the $k$ variables for bits of $i$. Note that $C_2$ can be encoded
by $\le m_{t-2}$ bits and $\tru_{s,k}(C_2) = b_2$.

Continuing in an analogous way for $t$ steps we show that the $t$-tuple of
functions (\ref{10.4.20a}) cannot have the claimed property. Note that the final $C_t$ witnessing that
$b_t \in rng(\tru_{s,k})$ too can be encoded by $m_0$ bits  and hence all
circuits $C_i$ have size at most  $m_0 = m^\epsilon = 2^{\epsilon k}$.

\qed

The reader who is confident that hypothesis (H) is false can interpret the theorem as saying that in order 
to disprove (H) it suffices to prove in $\tpv$ the existence of (a table of) a Boolean function
with an exponential circuit complexity.

\medskip

It would be desirable to prove the theorem for some other p-time function $g$
under a weaker or different hypothesis than (H). 
A good candidate for $g$ may be the proof complexity generator defined in \cite[Sec.3]{Kra-generator} (or see
\cite[Sec.19.4]{prf}). It is not difficult to modify the proof of Theorem \ref{8.4.20a} for this function (with suitable
parameters).

\end{document}